\documentclass[preprint,1p]{elsarticle}
\usepackage{graphics}
\usepackage{epsfig}
\usepackage{mathptmx}
\usepackage{amsmath}
\usepackage{amssymb}
\usepackage{hyperref}
\usepackage{fullpage}
\usepackage{amsthm}
\usepackage{lineno}

\journal{Name of Journal}

\begin{document}
\begin{frontmatter}
\title{Stancu-variant of generalized Baskakov operators}
\author[]{Nadeem Rao\corref{cor1}}
\ead{nadeemrao1990@gmail.com}
\author{Abdul Wafi}
\ead{abdulwafi2k2@gmail.com}
\address{\textbf{Department of Mathematics, Jamia Millia Islamia, New Delhi-110025, India}}
\cortext[cor1]{Corresponding Author}

\begin{abstract}
In the present paper, we introduce  Stancu-variant of generalized Baskakov operators and study the rate of convergence using modulus of continuity, order of approximation for the derivative of function $f$. Direct estimate is proved using K-functional and Ditzian-Totik modulus of smoothness. In the last, we have proved Voronovskaya type theorem.
\end{abstract}
\begin{keyword}
Generalized Baskakov operators, Modulus of continuity, Ditzian-Totik modulus of smoothness, Voronovskaya.
\end{keyword}
\end{frontmatter}
\section{Introduction}
In 1998 V. Mihesan\cite{Mihesa} constructed an important generalization of the well known Baskakov operators on $[0,\infty)$ with non-negative constant $a$ independent of $n$\\
\begin{equation}
  B_n^a(f;x)=\sum_{k=0}^{\infty}W_{n,k}^a(x)f\bigg(\frac{k}{n}\bigg),
\end{equation}
where $f\in C[0,\infty)$ and
 \begin{eqnarray}
 W_{n,k}^a(x)=e^{-\frac{ax}{1+x}}\frac{p_k(n,a)}{k!} \frac{x^k}{(1+x)^{k+n}},
 \end{eqnarray}
 such that $\sum\limits_{k=0}^{\infty}W_{n,k}^a(x)=1$ and $p_k(n,a)=\sum\limits_{i=0}^{\infty}\big( ^n _i\big)(n)_ia^{k-i},$ with $(n)_0=1,
 (n)_i=n(n+1)...(n+i-1).$\\ \\
In the last decade, many papers were published for generalized Baskakov operators on order of approximation, Voronovskaya type theorem, Kantorovich form, and order of approximation for the derivative of the function(\cite{Wafi1},\cite{wafi2}). \\
For $f\in C[0,\infty)$, Stancu\cite{D.D.} introduced the sequence of positive linear operators
\begin{eqnarray*}
  B_n^{\alpha,\beta}(f;x)=\sum_{k=0}^{n}p_{n,k}(x)f\bigg(\frac{k+\alpha}{n+\beta}\bigg),
\end{eqnarray*}
where $p_{n,k}(x)=(^n_i)x^k(1-x)^{n-k}$ and $\alpha,\beta$ are any two two non negative real numbers such that $0\leq\alpha\leq\beta$. If $\alpha=\beta=0$, it reduces to so called Bernstein operators. Recently, many researchers (\cite{ali},\cite{gupta1},\cite{gupta2},\cite{gupta3},\cite{radu},\cite{sucu},\cite{wang})have introduced the Stancu-variant for different linear positive operators. Motivated by the above development, we are giving a Stancu-variant of the operators (1) as\\
\begin{eqnarray}
L_{n,a}^{\alpha,\beta}(f;x)= \sum_{k=0}^{\infty}W_{n,k}^a(x)f\bigg(\frac{k+\alpha}{n+\beta}\bigg),
\end{eqnarray}
where $W_{n,k}^a(x)$ defined in (2) and $0\leq\alpha\leq\beta$. For $\alpha=\beta=0$, we get the operators (1). \section{\textrm{\textbf{ Approximation properties of $L_{n,a}^{\alpha,\beta}$}} }
To prove the approximation properties of $L_{n,a}^{\alpha,\beta}$, we need the following lemma\cite{Wafi1}\\ \\
\textbf{Lemma2.1} For $a,x\geq 0,n=1,2,...,$ we have\\
\begin{eqnarray*}
B_n^a(1;x)&=&1,\\
 B_n^a(t;x)&=& x+\frac{ax}{n(1+x)},\\
 B_n^a(t^2;x)&=& \frac{x^2}{n}+\frac{x}{n}+x^2+\frac{a^2x^2}{n^2(1+x)^2}+\frac{2ax^2}{n(1+x)}+\frac{ax}{n^2(1+x)},\\
 B_n^a(t^3;x)&=& x^3+\frac{3x^2(1+x)}{n}+\frac{x(1+x)(1+2x)}{n^2}+\frac{3ax^3}{n(1+x)}+\frac{1}{n^2}\bigg(3ax^2+\frac{3a^2x^3}{(1+x)^2}+\frac{3ax^2}{(1+x)}\bigg)\\
 &&+\frac{1}{n^3}\bigg(\frac{ax}{1+x}+\frac{3a^2x^2}{(1+x)^2}+\frac{a^3x^3}{(1+x)^3}\bigg),\\
 B_n^a(t^4;x)&=& x^4+\frac{6x^3(1+x)}{n}+\frac{x^2(1+x)(7+11x)}{n^2}+\frac{x(1+x)(6x^2+6x+1)}{n^3}+\frac{4ax^4}{n(1+x)}\\
 &&+\frac{1}{n^2}\bigg(\frac{12ax^4}{(1+x)+\frac{6a^2x^4}{(1+x)^2}}+\frac{18ax^3}{1+x}\bigg)+\frac{1}{n^3}\bigg(\frac{8ax^4}{1+x}+\frac{6a^2x^4}{(1+x)^2}+\frac{4a^3x^4}{(1+x)^3}+\frac{18ax^3}{(1+x)}+\frac{18a^2x^3}{(1+x)^2}+\frac{14ax^2}{(1+x)}\bigg)\\
 &&+\frac{1}{n^4}\bigg(\frac{ax}{1+x}+\frac{7a^2x^2}{(1+x)^2}+\frac{6a^3x^3}{(1+x)^3}+\frac{a^4x^4}{(1+x)^4}\bigg).
\end{eqnarray*}
Next, we prove\\ \\
\textbf{Lemma 2.2} Let $a,x\geq0$ and $n=1,2,3....$ Then for the operators defined in (3), we have
\begin{eqnarray*}
(i)\quad L_{n,a}^{\alpha,\beta}(1;x)&=&1,\\
(ii) \quad L_{n,a}^{\alpha,\beta}(t;x)&=&\frac{n}{n+\beta}x+\frac{a}{n+\beta}\frac{x}{1+x}+ \frac{\alpha}{n+\beta},\\
(iii)\quad L_{n,a}^{\alpha,\beta}(t^2;x)&=&\frac{n^2+n}{(n+\beta)^2}x^2+\frac{n(1+2\alpha)}{(n+\beta)^2}x+\frac{a^2}{(n+\beta)^2}\frac{x^2}{(1+x)^2}+\frac{2an}{(n+\beta)^2}\frac{x^2}{(1+x)}+\frac{a(1+2\alpha)}{(n+\beta)^2}\frac{x}{1+x} + \frac{\alpha^2}{(n+\beta)^2}\\
(iv) \quad L_{n,a}^{\alpha,\beta}(t^3;x)&=&\frac{n^3+3n^2+2n}{(n+\beta)^3}x^3+\frac{n^2(3+3\alpha)+n(3+3\alpha+3a)}{(n+\beta)^3}x^2+\frac{n(1+3\alpha+3\alpha^2)}{(n+\beta)^3}x+\frac{3an^2}{(n+\beta)^3}\frac{x^3}{(1+x)}\\
&&+\frac{n}{(n+\beta)^3}\bigg(\frac{3a^2x^3}{(1+x)^2}+\frac{3ax^2}{1+x}+\frac{6a \alpha x^2}{1+x} \bigg)+\frac{1}{(n+\beta)^3}\bigg(\frac{ax}{1+x}+\frac{3a^2x^2}{(1+x)^2}+\frac{a^3x^3}{(1+x)^3}\\
&&+\frac{3\alpha a^2x^2}{(1+x)^2}+\frac{3\alpha^2 ax}{1+x} +\alpha^3\bigg),\\
(v) \quad L_{n,a}^{\alpha,\beta}(t^4;x)&=&\frac{n^4+6n^3+11n^2+6n}{(n+\beta)^4}x^4+\frac{(6+4\alpha)n^3+(18+12\alpha)n^2+(9+8\alpha)n}{(n+\beta)^4}x^3+\bigg(\frac{(7+12\alpha+6\alpha^2)n^2}{(n+\beta)^4}\\
&&+\frac{(7+12\alpha+12\alpha a+6\alpha ^2)n}{(n+\beta)^4}\bigg)x^2+\frac{(1+4\alpha+6\alpha ^2+4\alpha ^3)n}{(n+\beta)^4}x+\frac{4an^3+12an^2+8an}{(n+\beta)^4}\frac{x^4}{1+x}\\
&&+\frac{6a^2n^2+6a^2n}{(n+\beta)^4}\frac{x^4}{(1+x)^2}+\frac{4a^3n}{(n+\beta)^4}\frac{x^4}{(1+x)^3}+\frac{a^4}{(b+\beta)^4}\frac{x^4}{(1+x)^4}+\frac{18an^2+18an}{(n+\beta)^4}\frac{x^3}{(1+x)}\\
&&+\frac{(18a^2+12a^2\alpha)n}{(n+\beta)^4}\frac{x^3}{(1+x)^2}+\frac{6a^3+4\alpha a^3}{(n+\beta)^4}\frac{x^3}{(1+x)^3}+\frac{(12a\alpha ^2+12a\alpha+14a)n}{(n+\beta)^4}\frac{x^2}{(1+x)}\\
&&+\frac{7a^2+12a^2\alpha+6a^2\alpha^2}{(n+\beta)^4}\frac{x^2}{(1+x)^2}+\frac{a+4\alpha a+6\alpha^2a+4\alpha^3a}{(n+\beta)^4}\frac{x}{1+x}+\frac{\alpha^4}{(n+\beta)^4}.
\end{eqnarray*}
\textbf{Proof} To prove these identities, we use the lemma(2.1) and linearity property
\begin{eqnarray*}
  L_{n,a}^{\alpha,\beta}(t;x)&=&\frac{n}{n+\beta}B_n^a(t;x)+\frac{\alpha}{n+\beta}B_n^a(1;x).\\
\end{eqnarray*}
In similar manner, we can prove identities (iii), (iv) and (v).\\ \\
\textbf{Lemma 2.3} Let $\psi_x^i(t)=(t-x)^i,i=1,2,3,...$. For $a,x\geq 0$ and $n=1,2,3...$,
\begin{eqnarray*}
L_{n,a}^{\alpha,\beta}(\psi_x^0(t);x)&=&1,\\
L_{n,a}^{\alpha,\beta}(\psi_x^1(t);x)&=&\bigg(\frac{n}{n+\beta}-1\bigg)x+\frac{a}{n+\beta}\frac{x}{1+x}+ \frac{\alpha}{n+\beta},\\
L_{n,a}^{\alpha,\beta}(\psi_x^2(t);x)&=&\frac{n+\beta^2}{(n+\beta)^2}x^2+\frac{n-2\alpha\beta}{(n+\beta)^2}x+\frac{a^2}{(n+\beta)^2}\frac{x^2}{(1+x)^2}-\frac{2a\beta}{(n+\beta)^2}\frac{x^2}{(1+x)}+\frac{a(1+2\alpha)}{(n+\beta)^2}\frac{x}{1+x} + \frac{\alpha^2}{(n+\beta)^2},\\
L_{n,a}^{\alpha,\beta}(\psi_x^4(t);x)&=&\frac{(3-12\beta)n^2+(6+4\beta+2\beta^2+4\beta^3)n+\beta^4}{(n+\beta)^4}x^4+\bigg(\frac{(6-12a-12\beta)n^2+(9+8\alpha-12\beta(1+a+\alpha+\alpha\beta))n}{(n+\beta)^4}\\
&&+\frac{(6-12a-12\beta-12\alpha\beta^2)}{(n+\beta)^4}\bigg)x^3+\frac{3n^2+(7-4\beta+12\alpha a-12\alpha\beta+6\alpha^2)n+6\alpha^2\beta^2}{(n+\beta)^4}x^2\\
&&+\frac{(1+4\alpha+6\alpha^2)n-4\alpha^3\beta}{(n+\beta)^4}x+\frac{12an^2+8an-4a\beta^3}{(n+\beta)^4}\frac{x^4}{(1+x)}+\frac{6a^2n+6a^2\beta^2}{(n+\beta)^4}\frac{x^4}{(1+x)^2} -\frac{4a^3\beta}{(n+\beta)^4}\frac{x^4}{(1+x)^3}\\
&&+\frac{a^4}{(n+\beta)^4}\frac{x^4}{(1+x)^4}+\frac{12an^2+18an+6a(1+2\alpha)\beta^2}{(n+\beta)^4}\frac{x^3}{1+x}+\frac{6a^2n-(12a^2+12\alpha a^2)\beta}{(n+\beta)^4}\frac{x^3}{(1+x)^2}\\
&&+\frac{(6a^3+4\alpha a^3)}{(n+\beta)^4}\frac{x^3}{(1+x)^3}+\frac{(12a\alpha+8a-6a\alpha^2)n-(6a+18\alpha^2a)\beta}{(n+\beta)^4}\frac{x^2}{1+x}+\frac{7a^2+12a^2\alpha+6a^2\alpha^2}{(n+\beta)^4}\frac{x^2}{(1+x)^2}\\
&&+\frac{(a)+4\alpha a+6\alpha^2a+4\alpha^3a}{(n+\beta)^4}\frac{x}{1+x}+\frac{\alpha^4}{(n+\beta)^4}.\\
\end{eqnarray*}
\textbf{Proof} In view of  lemma(2.2) and using equalities,
 \begin{eqnarray*}
 L_{n,a}^{\alpha,\beta}(\psi_x(t);x)&=&L_{n,a}^{\alpha,\beta}(t;x)-xL_{n,a}^{\alpha,\beta}(1;x),\\
 L_{n,a}^{\alpha,\beta}(\psi_x^2(t);x)&=&L_{n,a}^{\alpha,\beta}(t^2;x)-2xL_{n,a}^{\alpha,\beta}(t;x)+x^2L_{n,a}^{\alpha,\beta}(1;x),\\
 L_{n,a}^{\alpha,\beta}(\psi_x^4(t);x)&=&L_{n,a}^{\alpha,\beta}(t^4;x)-4xL_{n,a}^{\alpha,\beta}(t^3;x)+6x^2L_{n,a}^{\alpha,\beta}(t^2;x)+4x^3L_{n,a}^{\alpha,\beta}(t;x)+x^4L_{n,a}^{\alpha,\beta}(1;x).
 \end{eqnarray*}
 we get the proof of this lemma.\\
\textbf{Lemma 2.4} Let $\psi_x^i(t)=(t-x)^i,i=1,2,3,...$. For $a,x\geq 0$ and $n=1,2,3...$, we have
\begin{eqnarray*}
\lim\limits_{n \rightarrow \infty}nL_{n,a}^{\alpha,\beta}(\psi_x^1(t);x)&=&\alpha-\beta x+a\frac{x}{1+x},\\
\lim\limits_{n\rightarrow \infty}nL_{n,a}^{\alpha,\beta}(\psi_x^2(t);x)&=&x^2+x,\\
\lim\limits_{n\rightarrow\infty}n^2L_{n,a}^{\alpha,\beta}(\psi_x^4(t);x)&=&(3-12\beta)x^4+(6-12a-12\beta)x^3+3x^2+12a\frac{x^2}{1+x}+12a\frac{x^3}{1+x}.
\end{eqnarray*}
\section{The degree of approximation}
\textbf{Theorem 3.1} If $f\in C[0,\infty)$, $x\in[0,\infty)$ and $\omega(f;\delta)$ is the modulus of continuity, then\\
\begin{eqnarray*}
| L_{n,a}^{\alpha,\beta}(f;x)-f(x)|\leq\bigg\{1+\sqrt{\gamma_{n,a}^{\alpha,\beta}(x)}\bigg\}\omega(f;\delta_{n,\beta}),\\
    \end{eqnarray*}
where $\delta_{n,\beta}=(n+\beta)^{-\frac{1}{2}}$ and \\
\begin{eqnarray*}
\gamma_n^{\alpha,\beta}(x)=\frac{n+\beta^2}{n+\beta}x^2+\frac{n-2\alpha\beta}{n+\beta}x+\frac{a^2}{n+\beta}\frac{x^2}{(1+x)^2}-\frac{2a\beta}{n+\beta}\frac{x^2}{(1+x)}+\frac{a(1+2\alpha)}{n+\beta}\frac{x}{1+x} + \frac{\alpha^2}{n+\beta}.\\
\end{eqnarray*}
\textbf{Proof} Let $f\in C[0,\infty)$ and $x \geq 0$. Then, using linearity property and monotonicity of the operators defined by (3), we can easily find, for every $\delta> 0$, and $n\in N$,that
\begin{eqnarray*}
| L_{n,a}^{\alpha,\beta}(f;x)-f(x)| &\leq& \bigg\{1+\delta_{n,\beta}^{-1}\sqrt{L_{n,a}^{\alpha,\beta}(\psi_x^2;x)} \bigg\}\omega(f;\delta_{n,\beta}).\\
&\leq &\bigg \{ 1+\sqrt{ \frac{n+\beta^2}{n+\beta}x^2+\frac{n-2\alpha\beta}{n+\beta}x+\frac{a^2}{n+\beta}\frac{x^2}{(1+x)^2}-\frac{2a\beta}{n+\beta}\frac{x^2}{(1+x)}+\frac{a(1+2\alpha)}{n+\beta}\frac{x}{1+x} + \frac{\alpha^2}{n+\beta}}\bigg\}\omega(f;\delta_{n,\beta}),\\
\end{eqnarray*}
  which obtained by using Lemma 2.2 and choosing $\delta_{n,\beta}=(n+\beta)^{-\frac{1}{2}}$. Thus, we arrive at the result. \\ \\
\textbf{Remark} If we put $\alpha=\beta=0$, we find the same result given by Mihesan\cite{Mihesa}
\begin{eqnarray*}
| B_n^a(f;x)-f(x)| &\leq& \bigg\{1+\sqrt{x(1+x)+\frac{ax}{n(1+x)}\frac{(a+1)x+1}{(1+x)}} \bigg\}\omega(f;\delta),\\
\end{eqnarray*}
where $\delta=\frac{1}{\sqrt{n}}$, which shows that $ \delta_{n,\beta}\leq \delta $. Therefore, rate of convergence of $L_{n,a}^{\alpha,\beta}$ is better than $B_n^a$.\\ \\
Now, we will find the rate of convergence of operators defined by (3) in terms of modulus of continuity of first derivative of function $i.e.$  $\omega(f';\delta_{n,\beta})=\omega_1(f;\delta_{n,\beta})$, which is an improvement over the Theorem 3.1. This type of result was given for Bernstein polynomials by Lorentz (\cite{lorentz},p.p. 21).\\ \\
\textbf{Theorem 3.2} Let $f'(x)$ is the continuous derivative over $[0,\infty)$ and $\omega_1(f;\delta_{n,\beta})$ is the modulus of continuity of $f'(x)$. Then, for $a,x\geq 0$, $0\leq\alpha\leq\beta$, we have\\
\begin{eqnarray*}
 | L_{n,a}^{\alpha,\beta}(f;x)-f(x) | \leq\omega_1\big((n+\beta)^{-1}\big)\sqrt{L_{n,a}^{\alpha,\beta}(\psi_x^2(t);x)}\bigg\{1+\sqrt{(n+\beta)}\sqrt{L_{n,a}^{\alpha,\beta}(\psi_x^2(t);x)} \bigg\}.
\end{eqnarray*}
\textbf{Proof} For $x_1,x_2\in[a,b]$, we have\\
\begin{eqnarray}
\nonumber f(x_1)-f(x_2)&=&(x_1-x_2)f'(\xi),\\
 &=&(x_1-x_2)f'(x_1)+(x_1-x_2)[f'(\xi)-f'(x_1)],
\end{eqnarray}
where $x_1<\xi<x_2$. As we know that
\begin{eqnarray*}
 |(x_1-x_2)[f'(\xi)-f'(x_1)]| \leq| x_1-x_2 |(\lambda+1)\omega_1(\delta), \hspace{2cm} \lambda=\lambda(x_1,x_2;\delta).
\end{eqnarray*}
Next, we get
\begin{eqnarray}
| L_{n,a}^{\alpha,\beta}(f;x)-f(x) |= \bigg| \sum\limits_{n=0}^{\infty}W_{n,k}^a(x)\bigg\{f\bigg(\frac{k+\alpha}{n+\beta}\bigg)-f(x)\bigg\}\bigg |.
\end{eqnarray}
From (4) and (5), we obtained
\begin{eqnarray*}
| L_{n,a}^{\alpha,\beta}(f;x)-f(x) | &\leq& \bigg | \sum\limits_{n=0}^{\infty}W_{n,k}^{a}(x)\bigg(\frac{k+\alpha}{n+\beta}-x\bigg)f'(x) \bigg |+\omega_1(\delta_{n,\beta})\sum\limits_{k=0}^{\infty}\bigg | \frac{k+\alpha}{n+\beta}-x \bigg |(\lambda+1)W_{n,k}^{a}(x),\\
&\leq&\omega_1(\delta_{n,\beta})\bigg\{\sum\limits_{k=0}^{\infty}\bigg | \frac{k+\alpha}{n+\beta}-x \bigg | W_{n,k}^{a}(x)+\sum\limits_{\lambda\geq1}| \frac{k+\alpha}{n+\beta}-x| \lambda\bigg(x_1,\frac{k+\alpha}{n+\beta};\delta\bigg)W_{n,k}^{a}(x)\bigg\}\\
&\leq&\omega_1(\delta_{n\beta})\bigg\{\sum\limits_{k=0}^{\infty}\bigg | \frac{k+\alpha}{n+\beta}-x \bigg | W_{n,k}^{a}(x)+\delta^{-1}\sum\limits_{k=0}^{\infty}\bigg( \frac{k+\alpha}{n+\beta}-x\bigg)^2W_{n,k}^{a}(x)\bigg\}\\
&\leq&\omega_1(\delta_{n,\beta})\sqrt{L_{n,a}^{\alpha,\beta}(\psi_x^2(t);x)}\bigg\{1+\delta_{n,\beta}^{-1}\sqrt{L_{n,a}^{\alpha,\beta}(\psi_x^2(t);x)} \bigg\}.
\end{eqnarray*}
Taking $\delta_{n,\beta}=(n+\beta)^{-1}$, we get
\begin{eqnarray*}
 | L_{n,a}^{\alpha,\beta}(f;x)-f(x) | \leq\omega_1((n+\beta)^{-1})\sqrt{L_{n,a}^{\alpha,\beta}(\psi_x^2(t);x)}\bigg\{1+\sqrt{(n+\beta)}\sqrt{L_{n,a}^{\alpha,\beta}(\psi_x^2(t);x)} \bigg\},
\end{eqnarray*}
which is the required result.\\ \\
\section{Direct Estimate}
Here we introduced the Ditzian-Totik Modulus of smoothness\cite{dit1} which is defined as:
\begin{eqnarray*}
\omega^2_{\varphi^\lambda}(f;\delta)&=& \sup\limits_{0<h\leq \delta}\parallel \Delta^2_{h\varphi(x)}f(x) \parallel ,\\
&=& \sup\limits_{0<h\leq \delta} \hspace{0.5cm} \sup\limits_{x\pm h \varphi^\lambda \in [0,\infty)} |f(x-h\varphi^\lambda (x))-2f(x)+f(x+h\varphi^\lambda (x))|,\\
\end{eqnarray*}
where $\varphi^2(x)=x(1-x)$.
And, Peetre's K-functional is given by
\begin{eqnarray}
K_{\varphi^\lambda}(f,\delta^2)=\inf\limits_g \bigg(\|f-g\|_{C[0,\infty)} +\delta^2\|\varphi^2\lambda g''\|_{C[0,\infty)}\bigg), \hspace{2cm} g,g'\in AC_{loc}.
\end{eqnarray}
The K-functional is equivalent to the modulus of smoothness, i.e.,
\begin{eqnarray}
C^{-1}K_{\varphi^\lambda}(f,\delta^2)\leq \omega^2_{\varphi^\lambda}(f,\delta)\leq C K_{\varphi^\lambda}(f,\delta^2).
\end{eqnarray}
First result based on Ditziaz-Totik modulus of smoothness was given by Ditzian\cite{dit} for the Bernstein polynomials as:
\begin{eqnarray*}
| B_n(f;x)-f(x)| \leq C \omega^2_{\varphi^\lambda}(f,n^{-\frac{1}{2}} \varphi(x)^{1-\lambda}).
\end{eqnarray*}
Now, we prove the similar result for the operator $L_{n,a}^{\alpha,\beta}$. \\ \\ \\
\textbf{Theorem 4.1} For $a,x\geq 0$, and $0\leq \alpha\leq \beta$, we have\\
\begin{eqnarray*}
 | L_{n,a}^{\alpha,\beta}(f;x)-f(x)| \leq C \omega^2_{\varphi^\lambda}\big(f,(n+\beta)^{-\frac{1}{2}}  \varphi(x)^{1-\lambda}\big)\hspace{0.2cm} for\hspace{0.2cm} large \hspace{0.2cm}n.\\
 \end{eqnarray*}
 \textbf{Proof} Using (6),(7),we can choose $g_n\equiv g_{n,x,\lambda}$ for fixed $x$ and $\lambda+1$ such that
 \begin{eqnarray}
 \parallel\ f-g \parallel_{C[0,\infty)} \leq A\omega^2_{\varphi^\lambda}\big(f,n^{-\frac{1}{2}}\varphi(x)^{1-\lambda}\big),\\
 n^{-1}\varphi(x)^{2-2\lambda}\| \varphi^{2\lambda}g''\|_{C[0,\infty)}\leq B\omega^2_{\varphi^\lambda}\big(f,n^{-\frac{1}{2}}\varphi(x)^{1-\lambda}\big).
\end{eqnarray}
Next
\begin{eqnarray*}
 | L_{n,a}^{\alpha,\beta}(f;x)-f(x)| & \leq & | L_{n,a}^{\alpha,\beta}(f-g_n;x)-(f-g_n)(x) | +| L_{n,a}^{\alpha,\beta}(g_n;x)-g_n(x)|,\\ \\
 &\leq & 2\parallel f-g_n \parallel_{C[0,\infty)}+| L_{n,a}^{\alpha,\beta}(g_n;x)-g_n(x)|.
\end{eqnarray*}
From (8), we get
\begin{eqnarray}
 | L_{n,a}^{\alpha,\beta}(f;x)-f(x)| & \leq 2A\omega^2_{\varphi^\lambda}\big(f,n^{-\frac{1}{2}}\varphi(x)^{1-\lambda}\big) +| L_{n,a}^{\alpha,\beta}(g_n;x)-g_n(x)|.
\end{eqnarray}
Now, the last term can be calculated by using Taylor's formula
 \begin{eqnarray*}
 | L_{n,a}^{\alpha,\beta}(g_n(t)-g_n(x);x)| &\leq& | g'_n(x)L_{n,a}^{\alpha,\beta}((t-x);x)|+\Big| L_{n,a}^{\alpha,\beta}\bigg(\int\limits_{t}^{x}(x-u)g''_n(u)du;x\bigg)\Big|\\
 &\leq& L_{n,a}^{\alpha,\beta}\bigg( \frac{| x-\frac{k}{n}|}{\varphi^{2\lambda}(x)}\int\limits_{\frac{k}{n}}^x \varphi^{2\lambda}(u)|g''_n(u)du|;x \bigg)\\
 &\leq & \parallel \varphi^{2\lambda}g''_n\parallel_{C[0,\infty)}\frac{1}{\varphi^{2\lambda}(x)}L_{n,a}^{\alpha,\beta}((t-x)^2;x)\\
  &\leq &  \parallel \varphi^{2\lambda}g''_n\parallel_{C[0,\infty)}\frac{1}{\varphi^{2\lambda}(x)}\bigg[\frac{n+\beta^2}{(n+\beta)^2}x^2+\frac{n-2\alpha\beta}{(n+\beta)^2}x+\frac{a^2}{(n+\beta)^2}\frac{x^2}{(1+x)^2}-\frac{2a\beta}{(n+\beta)^2}\frac{x^2}{(1+x)}\\
 && +\frac{a(1+2\alpha)}{(n+\beta)^2}\frac{x}{1+x} + \frac{\alpha^2}{(n+\beta)^2}\bigg]\\
 &\leq &  \parallel \varphi^{2\lambda}g''_n\parallel_{C[0,\infty)}\frac{x(1+x)(n+\beta)^{-1}}{\varphi^{2\lambda}(x)}\bigg[\frac{n+\beta^2}{(n+\beta)}\frac{x}{1+x}
 +\frac{n-2\alpha\beta}{(n+\beta)}\frac{1}{1+x}+\frac{a^2}{(n+\beta)}\frac{x}{(1+x)^3}\\
 &&-\frac{2a\beta}{(n+\beta)}\frac{x}{(1+x)^2} +\frac{a(1+2\alpha)}{(n+\beta)}\frac{1}{(1+x)^2} + \frac{\alpha^2}{(n+\beta)x(1+x)}\bigg]\\
 &\leq&  \parallel \varphi^{2\lambda}g''_n\parallel_{C[0,\infty)}\varphi^{2-2\lambda}(x)(n+\beta)^{-1}\bigg[\frac{n+\beta^2}{(n+\beta)}\frac{x}{1+x}
 +\frac{n-2\alpha\beta}{(n+\beta)}\frac{1}{1+x}+\frac{a^2}{(n+\beta)}\frac{x}{(1+x)^3}\\
 &&-\frac{2a\beta}{(n+\beta)}\frac{x}{(1+x)^2} +\frac{a(1+2\alpha)}{(n+\beta)}\frac{1}{(1+x)^2} + \frac{\alpha^2}{(n+\beta)x(1+x)}\bigg].
 \end{eqnarray*}
 From (9), we have
 \begin{eqnarray}
| L_{n,a}^{\alpha,\beta}(g_n(t)-g_n(x);x)| &\leq& B\omega^2_{\varphi^\lambda}\big(f,(n+\beta)^{-\frac{1}{2}}\varphi(x)^{1-\lambda}\bigg[\frac{n+\beta^2}{(n+\beta)}\frac{x}{1+x}
 +\frac{n-2\alpha\beta}{(n+\beta)}\frac{1}{1+x}+\frac{a^2}{(n+\beta)}\frac{x}{(1+x)^3}\nonumber\\
&& -\frac{2a\beta}{(n+\beta)}\frac{x}{(1+x)^2} +\frac{a(1+2\alpha)}{(n+\beta)}\frac{1}{(1+x)^2} + \frac{\alpha^2}{(n+\beta)x(1+x)}\bigg].
 \end{eqnarray}
Using (10) and (11), we get\\
\begin{eqnarray*}
 | L_{n,a}^{\alpha,\beta}(f(t)-f(x);x)| &\leq& M\omega^2_{\lambda}\bigg(f,(n+\beta)^{\frac{-1}{2}}\varphi(x)^{1-\lambda}\bigg)\bigg[\frac{n+\beta^2}{(n+\beta)}\frac{x}{1+x}
 +\frac{n-2\alpha\beta}{(n+\beta)}\frac{1}{1+x}+\frac{a^2}{(n+\beta)}\frac{x}{(1+x)^3}\\
&& -\frac{2a\beta}{(n+\beta)}\frac{x}{(1+x)^2} +\frac{a(1+2\alpha)}{(n+\beta)}\frac{1}{(1+x)^2} + \frac{\alpha^2}{(n+\beta)x(1+x)}\bigg]\\
\end{eqnarray*}
where $M=max(2A,B)$. For a large value of $n$
\begin{eqnarray*}
 | L_{n,a}^{\alpha,\beta}(f(t)-f(x);x)| &\leq& M\omega^2_{\lambda}\bigg(f,(n+\beta)^{-\frac{1}{2}}\varphi(x)^{1-\lambda}\bigg).\\ \\
\end{eqnarray*}
Asymptotic relation is the study of rate of convergence for at least two times differentiable functions which was given by Voronovskaya \cite{voro}. Here, we prove a similar result.\\ \\
\textbf{\textbf{Theorem 4.2}} Let $a,x\geq0$, $ 0 \leq\alpha\leq\beta$ and $n\in N$. For $f \in C^2[0,\infty)$,  we have
 \begin{eqnarray*}
\lim\limits_{n\rightarrow\infty} n \{L_{n,a}^{\alpha,\beta}(f;x)-f(x)\}=\bigg(\alpha-\beta x+\frac{ax}{1+x}\bigg)f'(x)+\frac{x^2+x}{2}f''(x).
 \end{eqnarray*}
 \textbf{Proof} Let $x,t\in [0,\infty)$, $f\in C^2[0,\infty)$. By Taylor's formula, we have\\
 \begin{eqnarray*}
 f(t)=f(x)+(t-x)f'(x)+\frac{(t-x)^2}{2}f''(x)+\eta(t,x)(t-x)^2,
 \end{eqnarray*}
 where the function $\eta(t,x) \in C[0,\infty)$ and $\lim\limits_{t\rightarrow x} \eta(t,x)=0$. Multiplying both sides by $W_{n,k}^a(x)$ and summing over $k$, we get\\
 \begin{eqnarray*}
 L_{n,a}^{\alpha,\beta}(f;x)=f(x)L_{n,a}^{\alpha,\beta}(1;x)+f'(x)L_{n,a}^{\alpha,\beta}(t-x;x)+\frac{f''(x)}{2}L_{n,a}^{\alpha,\beta}((t-x)^2;x)+L_{n,a}^{\alpha,\beta}(\eta(t,x)(t-x);x).
 \end{eqnarray*}
 Using lemma(2.2), we obtain\\
 \begin{eqnarray}
\lim\limits_{n\rightarrow\infty} n \{L_{n,a}^{\alpha,\beta}(f;x)-f(x)\}=\bigg(\alpha-\beta x+\frac{ax}{1+x}\bigg)f'(x)+\frac{x^2+x}{2}f''(x)+\lim\limits_{n\rightarrow\infty}nL_{n,a}^{\alpha,\beta}(\eta(t;x)(t-x)^2;x).
 \end{eqnarray}
 Now, the last term can be obtained using Holder's inequality and lemma 2.4
 \begin{eqnarray*}
  nL_{n,a}^{\alpha,\beta}(\eta(t;x)(t-x)^2;x)\leq n^2L_{n,a}^{\alpha,\beta}((t-x)^4;x)L_{n,a}^{\alpha,\beta}(\eta(t;x)^2;x),\\
 \end{eqnarray*}
Let $\varphi(t;x)=\eta^2(t;x)$. Then, $\lim\limits_{t \rightarrow x}\varphi(t;x)=0$. Therefore,
\begin{eqnarray*}
  \lim\limits_{n\rightarrow\infty}nL_{n,a}^{\alpha,\beta}(\eta(t;x)(t-x)^2;x)=0.
\end{eqnarray*}
On substituting this value in equation (12) , we get the desired result.

\end{document}